\RequirePackage{fix-cm}
\documentclass{amsart}       
\usepackage{graphicx}
\usepackage{amsmath}
\usepackage{amssymb}
\usepackage{amsfonts}
\usepackage{latexsym}
\usepackage{oldgerm}

\newcommand{\smallmattwo}[4]{\left(\begin{smallmatrix} #1 & #2 \\ #3 &#4			       \end{smallmatrix}\right)}

\newcommand{\mattwo}[4]{\left(\begin{array}{cc} #1 & #2 \\ #3 &#4			       \end{array}\right)}

\newcommand{\vecttwo}[2]{\left(\begin{array}{c} #1 \\ #2
\end{array}\right)} 
\newcommand{\listtwo}[2]{\left\{\begin{array}{c} #1  \\
#2 \end{array}\right.}
\newcommand{\listtwotwo}[4]{\left\{\begin{array}{cc} #1 & #2 \\
#3 & #4 
\end{array}\right.}

\begin{document}

\title[Twisted Koecher-Maa{\ss} series of the Duke-Imamoglu-Ikeda lift]
 {Explicit formulas for the twisted Koecher-Maa{\ss} series
 of the Duke-Imamoglu-Ikeda lift and their applications}
\thanks{The author was partly supported by JSPS KAKENHI Grant Number 24540005, JSPS.}

\author{Hidenori Katsurada}


\maketitle

\begin{abstract}
We give an explicit formula for the twisted Koecher-Maa{\ss} series of the Duke-Imamoglu-Ikeda lift. As an application we 
prove a certain algebraicity result for the values of twisted Rankin-Selberg series at integers of half-integral weight modular forms,
which was not treated by Shimura \cite{Sh3}. 


\bigskip


\bigskip

\noindent
{\bf Mathmatical Subject Classification(2000)} 11F67
\end{abstract}

\section{Introduction}

It is an interesting problem to give an explicit formula for the Koecher-Maa{\ss} series of a Siegel modular form $F$ for 
the symplectic group $Sp_n({\bf Z})$, and 
several  results have been obtained (cf.   B\"ocherer \cite{B},  Ibukiyama and Katsurada  \cite{I-K1}, \cite{I-K2}, \cite{I-K3}).  Such explicit formulas are not only interesting in its own right but also have some important applications in the theory of modular forms. For example, we  refer to \cite{B-S},  \cite{DI}. 
 Now we consider a twist of such a Koecher-Maa{\ss}  series by a Dirichlet character $\chi.$
As for this,  in view of Saito \cite{Sai1} for example, we  can naturally consider the following Dirichlet series:
\[L^*(s,F,\chi)=\sum_T {\displaystyle  \chi(2^{2[n/2]}\det T)c_F(T) \over \displaystyle e(T) (\det T)^s},\] 
where $T$ runs over a complete set of representatives of $SL_{n}({\bf Z})$-equivalence classes of positive definite half-integral symmetric matrices of degree $n$, $c_F(T)$ is the $T$-th Fourier coefficient of $F$ and $e(T)=\#\{U \in SL_n({\bf Z});
 T[U]=T \}$. We  will sometimes call $L^*(s,F,\chi)$ the twisted Koecher-Maa{\ss} series of the second kind. 

On the other hand, Choie and Kohnen \cite{C-K} introduced a different type of ``twist". 
For a positive integer $N$, let $SL_{n,N}({\bf Z})=\{U \in SL_n({\bf Z}); U \equiv 1_n \ {\rm mod} \ N \}$ and $e_N(T)=\#\{U \in SL_{n,N}({\bf Z}); T[U]=T \}$. For a primitive Dirichlet character $\chi$ mod $N$,  the Koecher-Maa{\ss} series $L(s,F,\chi)$ of $F$ twisted by $\chi$ is defined to be   
\[\displaystyle L(s,F,\chi)=\sum_{T} {\displaystyle  \chi({\rm tr}(T))c_F(T) \over \displaystyle e_N(T) (\det T)^s},\]
where $T$ runs over a complete set of representatives
of $SL_{n,N}({\bf Z})$-equivalence classes of positive definite half-integral symmetric 
matrices of degree $n$. In \cite{C-K}, Choie and Kohnen proved an analytic continuation of $L(s,F,\chi)$ to the whole $s$-plane and a functional equation (cf. Theorem 2.1). 
Moreover they got a result on the algebraicity  of its special values (cf. Theorem 2.2.) We shall call $L(s,F,\chi)$ the twist of the first kind. 

In this paper we give explicit formulas for the twisted Koecher-Maa{\ss} series of the first and second kinds  associated with the Duke-Imamoglu-Ikeda lift and apply them to the study of the special values of the Rankin-Selberg series for  half-integral weight modular forms.  
We explain our main results more precisely. Let $k$ and $n$ be positive even integers such that $n \ge 4$ and $2k-n \ge 12.$ For  a cuspidal Hecke eigenform $h$  in the Kohnen plus subspace of weight $k-n/2+1/2$ for $\Gamma_0(4),$ let 
$I_n(h)$ be the Duke-Imamoglu-Ikeda lift of $h$ to the space of cusp forms of weight $k$ for $Sp_n({\bf Z}).$ Moreover let $S(h)$ be the normalized Hecke eigenform of weight  $2k-n$ for $SL_2({\bf Z})$ corresponding to $h$ under the 
Shimura correspondence, and $E_{n/2+1/2}$ be Cohen's Eisenstein series of weight $n/2+1/2$ for $\Gamma_0(4).$ 
We then give explicit formulas for $L(s,I_n(h),\chi)$ and $L^*(s,I_n(h),\chi)$  in terms of the twisted  Rankin-Selberg series  $R(s,h,E_{n/2+1/2},\eta)$ of $h$ and $E_{n/2+1/2}$ and twisted Hecke's $L$-function $L(s,S(h), \eta') $ of 
$S(h),$ where $\eta$ and $\eta'$ are Dirichlet characters related with $\chi.$ 
It is relatively  easy to get an explicit form of $L^*(s,I_n(h),\chi)$. In fact, by using the same argument as in Ibukiyama and Katsurada \cite{I-K2}, we can easily obtain its explicit formula (cf. Theorem 4.1). On the other hand, it seems nontrivial to get that of $L(s,I_n(h),\chi)$ (cf. Theorem 6.1), and we need some explicit formula for a certain character sum  associated with a Dirichlet character (cf. Theorem 5.6).  
Using Theorem 6.1 combined  with the result of Choie-Kohnen, we prove certain algebraicity results 
 on $R(s,h,E_{n/2+1/2},\eta)$ at an integer $s=m$ 
 (cf. Theorems 7.1 and 7.2), which were announced in \cite{Ka}.  We note that the algebraicity of the special values of such a Rankin-Selberg series at half-integers was investigated by Shimura \cite{Sh3}. However there are few results on the algebracity of such values at integers. As an attempt, Mizuno and the author \cite{K-M} proved linear dependency of Rankin-Selberg $L$-values of a cuspidal Hecke eigenform belonging to Kohnen plus subspace of half integral weight and the Zagier's Eisenstein series of weight $3/2.$ Our present result can be regarded as a generalization of our previous result.

{\bf Notation}  We denote by ${\bf e}(x)=\exp(2 \pi \sqrt{-1}x)$ for a complex number $x.$  For a
commutative ring $R$, we denote by $M_{mn}(R)$ the set of
$(m,n)$-matrices with entries in $R.$ In particular put $M_n(R)=M_{nn}(R).$   For an $(m,n)$-matrix $X$ and an $(m,m)$-matrix $A$, we write $A[X] = {}^t \! X A X,$ where $^t \!X$ denotes the
transpose of $X$.  Let $a$ be an element of $R.$ Then
for an element $X$ of $M_{mn}(R)$ we often use the same
symbol $X$ to denote the coset $X \ {\rm mod} \ aM_{mn}(R).$ Put $GL_m(R)= \{A \in M_m(R) \ | \ \det A  \in R^* \},$  and $SL_m(R) = \{A \in M_m(R) \ | \ \det A =1 \},$ where $\det
A$ denotes the determinant of a square matrix $A$ and $R^*$ is the unit group of $R.$ We denote by $S_n(R)$ the set of symmetric matrices of degree $n$ with entries in $R.$   
For a subset $S$ of $M_n(R)$ we denote by $S^{\times}$ the subset of $S$ consisting of non-degenerate matrices.  In particular, if
$S$ is a subset of $S_n({\bf R})$ with ${\bf R}$ the field of real numbers, we denote by $S_{>0}$ (resp. $S_{\ge 0}$) the subset of $S$
consisting of positive definite (resp. semi-positive definite) matrices. The group $SL_n({\bf Z})$ acts on the set $S_n({\bf R})$ in the following way:
  $$ SL_n({\bf Z}) \times S_n({\bf R}) \ni (g,A) \longrightarrow {}^tg Ag \in S_n({\bf R}).$$
  Let $G$ be a subgroup of $GL_n({\bf Z}).$ For a subset ${\mathcal B}$ of $S_n({\bf R})$ stable under the action of $G$ we denote by ${\mathcal B}/G$ the set of equivalence classes of ${\mathcal B}$ with respect to $G.$ We sometimes identify ${\mathcal B}/G$ with a complete set of representatives of ${\mathcal B}/G.$ Two symmetric matrices $A$ and $A'$ with
entries in $R$ are said to be equivalent   with respect to $G$ and write $A \sim_G A'$ if there is an element $X$ of $G$ such that $A'=A[X].$ 
For an integral domain $R$ of charactersitic different from $2$, let ${\mathcal L}_n(R)$
denote the set of half-integral matrices of degree $n$ over
$R,$ that is, ${\mathcal L}_n(R)$ is the set of symmetric
matrices of degree $n$ with entries in the field of fractions of $R$  whose $(i,j)$-component belongs to
$R$ or ${1 \over 2}R$ according as $i=j$ or not. In particular we put ${\mathcal L}_n={\mathcal L}_n({\bf Z}).$
For square matrices $X$ and $Y$ we write $X \bot Y = \mattwo{X}{O}{O}{Y}.$ For a subset $S$ of a ring $R$ we put  
$S^{\Box} =\{ s^2 \ | \ s \in S \}.$ 

\section{Twisted Koecher-Maa{\ss} series}
Put $J_n=\mattwo{O_n}{-1_n}{1_n}{O_n},$ where $1_n$ and $O_n$ denotes the unit matrix and the zero matrix of degree $n$, respectively.  Furthermore, put 
$$Sp_n({\bf Z})=\{M \in GL_{2n}({\bf Z})   \ | \  J_n[M]=J_n \}. $$
 Let $l$ be an integer or a half-integer, and $N$ a positive integer. Let $\Gamma_0^{(n)}(N)$ be the congruence subgroup of $Sp_n({\bf Z})$ consisting of matrices whose left lower $n \times n$ block are congruent to $O_n$ mod $N.$ Moreover let 
  $\chi$ be a Dirichlet character mod $N.$  We  then denote by ${\textfrak M}_l(\Gamma_0^{(n)}(N),\chi)$ the space of modular forms of weight $l$ and character $\chi$ for $\Gamma_0^{(n)}(N),$ and by ${\textfrak S}_l(\Gamma_0^{(n)}(N),\chi)$ the subspace of ${\textfrak M}_l(\Gamma_0^{(n)}(N),\chi)$ consisting of cusp forms. If $\chi$ is the trivial character mod $N$, we simply write ${\textfrak M}_l(\Gamma_0^{(n)}(N),\chi)$ and ${\textfrak S}_l(\Gamma_0^{(n)}(N),\chi)$ as ${\textfrak M}_l(\Gamma_0^{(n)}(N))$ and ${\textfrak S}_l(\Gamma_0^{(n)}(N)),$ respectively. Let $k$ be a positive integer, and let 
$F(Z) \in  {\textfrak M}_k(Sp_n({\bf Z})).$ Then $F(Z)$ has the Fourier expansion:
$$F(Z)= \sum_{T \in{{\mathcal  L}_n}_{\ge  0}} c_F(T) {\bf e}({\rm tr}(TZ)),$$
where ${\rm tr}(X)$ denotes the trace of a matrix $X.$ For $N \in {\bf Z}_{>0},$ put
$SL_{n,N}({\bf Z})=\{U \in SL_n({\bf Z}) \ | \ U \equiv 1_n \ {\rm mod} \ N \},$ and
for $T \in {{\mathcal L}_n}_{>0}$ put $e_N(T)=\#\{U \in SL_{n,N}({\bf Z}) \ | \ T[U]=T \}.$
For  a primitive Dirichlet character $\chi$ mod $N$ let 
$$L(s,F,\chi)=\sum_{T \in {{\mathcal L}_n}_{>0}/SL_{n,N}({\bf Z})} { \chi({\rm tr}(T))c_F(T) \over e_N(T) (\det T)^s}$$  
be the twisted Koecher-Maa{\ss} series of the first kind  of $F$ as in Section 1.  The following two theorems are due to Choie and Kohnen \cite{C-K}.

\bigskip

\noindent
{\bf Theorem 2.1.} {\it 
Let $F \in {\textfrak S}_k(Sp_n({\bf Z})),$ and $\chi$ a primitive character of conductor $N.$  Put
$$\gamma_n(s)=(2\pi)^{-ns}\prod_{i=1}^n \pi^{(i-1)/2}{\rm \Gamma}(s-(i-1)/2),$$
and 
$$\Lambda(s,F,\chi)=N^{ns}\tau(\chi)^{-1} \gamma_n(s)L(s,F,\chi) \quad ({\rm Re}(s) >>0),$$
where $\tau(\chi)$ is the Gauss sum of $\chi,$ and ${\rm \Gamma}(s)$ is the Gamma function.
Then $\Lambda(s,F,\chi)$ has an analytic continuation to the whole $s$-plane and has
the following functional equation:
$$\Lambda(k-s,F,\chi)=(-1)^{nk/2}\chi(-1)\Lambda(s,F,\overline \chi).$$
}

{\bf Theorem 2.2.}   {\it
 Let $F$ and $\chi$ be as above. Then there exists a finite dimensional $\overline {{\bf Q}}$-vector space $V_F$ in ${\bf C}$ such that
$$ L(m,F,\chi) \pi^{-nm} \in V_F$$
for any primitive character $\chi$ and any integer $m$ such that  $(n+1)/2 \le m \le k-(n+1)/2.$
}

Now let 
$$L^*(s,F,\chi)=\sum_{T  \in {{\mathcal L}_n}_{>0}/SL_n({\bf Z})}{\displaystyle  \chi(2^{2[n/2]}\det T)c_F(T) \over \displaystyle e(T) (\det T)^s}$$ 
be the twisted Koecher-Maa{\ss} series of the second kind of $F$ as in Section 1.  We will discuss a relation between these two Dirichlet series in Section 5.

\section{Review on the algebraicity of L-values of  elliptic modular forms of integral  and half-integral weight}

In this section, we review on the special values of L functions of elliptic modular forms of integral and half-integral weights. For a modular form $g(z)$ of integral or half-integral weight for a certain congruence subgroup $\Gamma$ of $SL_2({\bf Z}),$ 
 let ${\bf Q}(g)$ denote the field generated over ${\bf Q}$ by all the Fourier coefficients of $g,$ and for a Dirichlet character $\eta$ let ${\bf Q}(\eta)$ denote the field generated over ${\bf Q}$ by all the values of $\eta.$

First let 
$$f(z)=\displaystyle \sum_{m=1}^{\infty} c_f(m){\bf e}(mz)$$ be a normalized Hecke eigenform in ${\textfrak S}_{k}(SL_2({\bf Z})),$
and $\chi$ be a primitive Dirichlet character. Then let us define Hecke's $L$-function $L(s,f,\chi)$ of $f$ twisted by $\chi$ as 
$$L(s,f,\chi)=\displaystyle \sum_{m=1}^{\infty} c_f(m)\chi(m)m^{-s}.$$
Then we have the following result (cf. \cite {Sh2}):

\bigskip

\noindent
{\bf Proposition  3.1.} {\it  
There exist complex numbers   $u_{\pm}(f)$
 uniquely determined up to ${\bf Q}(f)^{\times}$ multiple such that $${L(m,f,\chi) \over (2\pi \sqrt{-1})^{m}\tau(\chi)u_{j}(f)} \in {\bf Q}(f){\bf  Q}(\chi)$$
for any integer $0< m \le k-1$ and a  primitive character $\chi,$
 where  $\tau(\chi)$ is the Gauss sum of $\chi$, and $j=+$ or $-$ according as $(-1)^{m}\chi(-1)=1$ or $-1.$ 
 }

\noindent
{\bf Corollary.} {\it Under the above notation and the assumption, we have
$$L(m,f,\chi) \pi^{-m}  \in \overline{{\bf Q}} u_{j}(f)$$
for any integer $0< m \le k-1$ and a  primitive character $\chi.$}
  
We remark that we have $L(m,f,\chi) \not=0$ if $m \not=k/2,$ and  $L(k/2,f,\chi) \not=0$ for infinitely many $\chi.$

Next let us consider the half-integral weight case. 
From now on we simply write $\Gamma_0^{(1)}(M)$ as $\Gamma_0(M).$ Let 
$$h_1(z)=\displaystyle \sum_{m=1}^{\infty} c_{h_1}(m){\bf e}(mz)$$ 
be a Hecke eigenform in ${\textfrak S}_{k_1+1/2}(\Gamma_0(4)),$ and 
$$h_2(z)=\displaystyle \sum_{m=0}^{\infty} c_{h_2}(m){\bf e}(mz)$$
be an element of  ${\textfrak M}_{k_2 +1/2}(\Gamma_0(4)).$ 
For a fundamental discriminant $D$ let $\chi_D$ be the Kronecker character corresponding to $D.$ 
Let $\chi$ be a primitive character mod $N.$ Then we define 
$$\widetilde R(s,h_1,h_2,\chi)=L(2s-k_1-k_2+1,\omega) \displaystyle \sum_{m=1}^{\infty} c_{h_1}(m)c_{h_2}(m) \chi(m)m^{-s},$$
where $\omega(d)=\chi_{-4}^{k_1-k_2}\chi^2(d).$
We also define $R(s,h_1,h_2,\chi)$ as 
 $$R(s,h_1,h_2,\chi)=L(2s-k_1-k_2+1,\chi^2) \displaystyle \sum_{m=1}^{\infty} c_{h_1}(m)c_{h_2}(m) \chi(m)m^{-s}.$$
Now let $S(h_1)$ be the normalized Hecke eigenform in ${\textfrak S}_{2k_1}(SL_2({\bf Z}))$ corresponding to $h_1$ under the Shimura correspondence. Then the following result is due to Shimura \cite {Sh3}.

\bigskip

\noindent
{\bf Proposition 3.2.} {\it 
 Assume that $k_1 > k_2.$  Under the above notation we have 
$${\widetilde R(m+1/2,h_1,h_2,\chi) \over u_{-}(S(h_1)) \tau(\chi^2)\pi^{-k_2+1+2m} \sqrt{-1}} \in {\bf Q}(h_1){\bf Q}(h_2){\bf Q}(\chi)$$ 
 for any  integer $k_2 \le m \le k_1-1$ and a primitive character $\chi.$
}
\begin{proof} Let $N$ be the conductor of $\chi.$ Put
$${h_2}_{\chi}(z)=\sum_{m=0}^{\infty} c_{h_2}(m)\chi(m){\bf e}(mz).$$
Then ${h_2}_{\chi}(z) \in {\textfrak M}_{k_2+1/2}(4N^2,\chi^2).$ We can regard $h_1$ as an element of ${\textfrak S}_{k_1+1/2}(\Gamma_0(4N^2)).$ Then the assertion follows from [\cite{Sh3}, Theorem 2].
\end{proof}
\bigskip

\noindent
{\bf Corollary.} {\it Assume that $c_{h_1}(n), c_{h_2}(n) \in \overline { {\bf Q}}$ for any $n \in {\bf Z}_{\ge 0}.$ Then there 
exists a one-dimensional  $\overline {{\bf Q}}$-vector space $U_{h_1,h_2}$ in $ {\bf C}$ such that 
$$\widetilde R(m+1/2,h_1,h_2,\chi) \pi^{-2m} \in U_{h_1,h_2}$$ 
for any integer $k_2 \le m \le k_1-1$ and a primitive character $\chi.$}

\bigskip

\section {Explicit formulas for the twisted Koecher-Maa{\ss}  series of the second kind of the  Duke-Imamoglu-Ikeda  lift}

Throughout this section, we assume that $n$ and $k$ are even positive integers.
Let  $h$ be a Hecke eigenform of weight $k-n/2+1/2$ for $\Gamma_0(4)$ belonging to the Kohnen plus space. Then $h$  has the following Fourier expansion:  
$$h(z)=\sum_{e}c_h(e){\bf e}(ez),$$ 
where $e$ runs over all positive integers such that $(-1)^{k-n/2}e \equiv 0, 1 \ {\rm mod} \ 4.$
Let 
$$S(h)(z)=\sum_{m=1}^{\infty}c_{S(h)}(m){\bf e}(mz)$$
be the  normalized Hecke eigenform of weight $2k-n$ for $SL_2({\bf Z})$
 corresponding to $h$ via  the Shimura correspondence (cf. \cite{Ko}.)   
For a prime number $p$ let $\beta_p$ be a nonzero complex number such that
$\beta_p+\beta_p^{-1}=p^{-k+n/2+1/2}c_{S(h)}(p).$
For non-negative integers $l$ and  $m$, the Cohen function $H(l,m)$ is given by $H(l,m)=L_{-m}(1-l)$. Here
\begin{eqnarray*}
&\,&L_{D}(s)\\&=&\left\{
\begin{array}{rl}
\displaystyle{\zeta(2s-1)},&\quad\mbox{$D=0$}\\
\displaystyle{L(s, \chi_{D_K})\sum_{a|f}\mu(a)\chi_{D_K}(a)a^{-s}\sigma_{1-2s}(f/a)},&\quad\mbox{$D\ne 0,$~$D\equiv 0, 1\ {\rm mod} \ 4$}\\
\displaystyle{0},&\quad\mbox{$D\equiv 2, 3\ {\rm mod} \ 4$,}
\end{array}\right.
\end{eqnarray*}
where the positive integer $f$ is defined by $D=D_{K}f^2$ with the discriminant $D_{K}$ of $K={\mathbf Q}(\sqrt{D})$, 
$\mu$ is the M$\ddot{\rm o}$bius function, and $\sigma_{s}(n)=\sum_{d|n}d^s$. 
Furthermore, for an even integer $l \ge 4,$  we define the Cohen Eisenstein series $E_{l+1/2}(z)$ 
 by
 $$E_{l+1/2}(z)=\sum_{e=0}^{\infty}H(l,e) {\bf e}( ez).$$
 It is known that  $E_{l+1/2}(z)$ is a modular form of weight $l+1/2$  for $\Gamma_0(4)$ belonging to the Kohnen plus space. 

For a prime number $p$ let
 ${\bf Q}_p$ and ${\bf Z}_p$ be the field of $p$-adic numbers, and the ring of $p$-adic integers, respectively. We denote by $\nu_p$ the additive valuation on ${\bf Q}_p$ normalized so that $\nu_p(p)=1,$ and by ${\bf e}_p$ the continuous homomorphism from the additive group ${\bf Q}_p$ to ${\bf C}^{\times}$ such that ${\bf e}_p(a)={\bf e}(a)$ for $a \in {\bf Q}.$  We put ${\mathcal L}_{n,p}={\mathcal L}_n({\bf Z}_p).$ We also put $S_n({\bf Z}_p)_e=2{\mathcal L}_{n,p}$ and $S_n({\bf Z}_p)_{o}=S_n({\bf Z}_p) \smallsetminus S_n({\bf Z}_p)_e.$  For a $p$-adic number $c$ put $$\widetilde \xi_p(c)=1,-1 \ {\rm
or} \ 0 $$ according as ${\bf Q}_p(\sqrt c)={\bf Q}_p,{\bf Q}_p(\sqrt c)/{\bf
Q}_p$ is quadratic unramified, or ${\bf Q}_p(\sqrt c)/{\bf Q}_p$ is quadratic
ramified.  We note that $\tilde \xi_p(D)=\chi_D(p)$ for a fundamental discriminant $D.$ For a non-degenerate half-integral matrix $T$ over ${\bf Z}_p,$ let 
$$b_p(T,s)=\sum_{R \in S_n({\bf Q}_p)/S_n({\bf Z}_p)} {\bf e}_p( {\rm tr}(TR))p^{-\nu_p(\mu_p(R))s}$$ 
 be the local Siegel series, where $\mu_p(R)=[R{\bf Z}_p^n+{\bf Z}_p^n:{\bf Z}_p^n].$  Then there exists a polynomial 
 $F_p(T,X)$ in $X$ such that
 $$b_p(T,s)=F_p(T,p^{-s})(1-p^{-s})(1-\xi_p(T)p^{n/2-s})^{-1}\prod_{i=1}^{n/2} 
(1-p^{2i-2s})$$ 
(cf. \cite{Ki1},)  where $\xi_p(T)= \widetilde \xi_p((-1)^{n/2} \det T).$  For a positive definite half integral matrix $T$ of degree $n$ write $(-1)^{n/2}\det (2T)$ as $(-1)^{n/2}\det (2T)
 = {\textfrak d}_T {\textfrak f}_T^2$ with ${\textfrak d}_T$ a fundamental discriminant and ${\textfrak f}_T$ a positive integer. 
We then put  
$$\widetilde F_p(T,X)=X^{-\nu_p({\textfrak f}_T)}F_p(T,p^{-(n+1)/2}X),$$
and 
 $$c_{I_n(h)}(T)=c_h(|{\textfrak d}_T|)  \prod_p (p^{k-n/2-1/2})^{\nu_p({\textfrak f}_T)}\widetilde F_p(T,\beta_p).$$
We note that $c_{I_n(h)}(T)$ does not depend on the choice of $\beta_p.$ 
Define a Fourier series $I_n(h)(Z)$ by 
$$I_n(h)(Z)= \sum_{T \in {{ \mathcal  L}_n}_{> 0}} c_{I_n(h)}(T){\bf e}({\rm tr}(TZ)).$$ 
In \cite{I} Ikeda showed that $I_n(h)(Z)$ is a  Hecke eigenform in ${\textfrak S}_k(Sp_n({\bf Z}))$ and its standard $L$-function  
$L(s,I_n(h),{\rm St})$ is given by 
$$L(s,I_n(h),{\rm St})=\zeta(s)\prod_{i=1}^nL(s+k-i,S(h)).$$
We call $I_n(h)$ the Duke-Imamoglu-Ikeda lift (D-I-I lift) of $h.$

\bigskip

\noindent
{\bf Theorem 4.1.} {\it 
Let $\chi$ be a primitive Dirichlet character mod $N$. Then we have 
$$L^*(s,I_n(h),\chi)= c_{n}R(s,h,E_{n/2+1/2}, \chi )\prod_{j=1}^{n/2-1} L(2s-2j,S(h),\chi^2) $$
$$+d_{n}c_h(1)  \prod_{j=1}^{n/2} L(2s-2j+1,S(h),\chi^2), $$
where $c_{n}$ and $d_{n}$ are nonzero rational numbers depending only on $n.$ 
 } 

\bigskip

To prove Theorem 4.1, we reduce the problem to local computations. For $a,b \in {\bf Q}_p^{\times}$ let $(a,b)_p$ the Hilbert symbol on ${\bf Q}_p.$ Following Kitaoka \cite{Ki2}, we define the Hasse invariant $\varepsilon(A)$ of $A \in S_m({\bf Q}_p)^{\times}$ by 
$$\varepsilon(A)=\prod_{1 \le i \le j \le n}(a_i,a_j)_p$$
if $A$ is equivalent to $a_1 \bot \cdots \bot a_n$ over ${\bf Q}_p$ with some $a_1,a_2,\cdots ,a_n \in {\bf Q}_p^{\times}.$ 
For  non-degenerate symmetric matrices $A$  of degree $n$ with entries in ${\bf Z}_p$  we define the local density $\alpha_p(A)=\alpha_p(A,A)$ 
representing $A$ by $A$ as
$$\alpha_p(A)=2^{-1}\lim_{a \rightarrow
\infty}p^{a(-n^2+n(n+1)/2)}\#{\mathcal A}_a(A,A),$$
where $${\mathcal A}_a(A,A)=\{X \in
M_{n}({\bf Z}_p)/p^aM_{n}({\bf Z}_p) \ | \ A[X]-B \in p^aS_n({\bf Z}_p)_e \},$$
Furthermore put 
$$M(A)=\sum_{A' \in {\mathcal G}(A)} {1 \over e(A')}$$
for a positive definite symmetric matrix $A$ of degree $n$ with entries in ${\bf Z},$ where ${\mathcal G}(A)$ denotes the set of $SL_{n}({\bf Z})$-equivalence classes belonging to the genus of $A.$   Then
by Siegel's main theorem on quadratic forms, we obtain 
$$M(A)=\kappa_{n} 2^{-n/2}\det A^{(n+1)/2} \prod_p \alpha_p(A)^{-1}$$
where $e_{n}=1$ or $2$ according as $n=1$ or not, and $\kappa_{n}={\rm \Gamma}_{\bf C}(n/2)\prod_{i=1}^{n/2-1}{\rm \Gamma}_{\bf C}(2i)$
with ${\rm \Gamma}_{\bf C}(s)=2(2\pi )^{-s}{\rm \Gamma}(s)$
(cf. Theorem 6.8.1 in \cite{Ki2}
).  Put  
$${\mathcal F}_{p}=\{d_0 \in {\bf Z}_p \ | \ \nu_p(d_0) \le 1\}$$ 
 if $p$ is an odd prime, and  
$${\mathcal F}_{2}=\{d_0 \in {\bf Z}_2 \ | \  d_0 \equiv 1 \ {\rm mod} \ 4, \ {\rm  or} \  d_0/4  \equiv -1 \  {\rm mod} \ 4,  \ {\rm or} \ \nu_2(d_0)=3 \}.$$
For $d \in {\bf Z}_p^{\times}$ put 
 \begin{eqnarray*}
\lefteqn{S_{n}({\bf Z}_p,d)}\\
=&\{T \in S_{n}({\bf Z}_p) \ |\ (-1)^{n/2} \det T=p^{2i}d \ {\rm mod} \ {{\bf Z}_p^*}^{\Box} \ {\rm with \ some \ } \ i \in {\bf Z} \}, 
\end{eqnarray*}
and $S_{n}({\bf Z}_p,d)_x=S_{n}({\bf Z}_p,d) \cap S_{n}({\bf Z}_p)_x$ for $x=e$ or $o.$   
 Put ${\mathcal L}_{n,p}^{(0)}=S_{n}({\bf Z}_p)_e^{\times}$
  and
${\mathcal L}_{n,p}^{(0)}(d)= S_{n}({\bf Z}_p,d) \cap {\mathcal L}_{n,p}^{(0)}.$   
Let $\iota_{n,p}$ be the constant function on ${\mathcal L}_{n,p}^{\times}$ taking the value 1, and $\varepsilon_{n,p}$ the function on ${\mathcal L}_{n,p}^{\times}$ assigning the Hasse invariant of $A$ for $A \in {\mathcal L}_{n,p}^{\times}.$  We sometimes drop the suffix and write $\iota_{n,p}$ as $\iota_p$ or $\iota$ and the others if there is no fear of confusion.
From now on we sometimes write $\omega=\varepsilon^l$ with $l=0$ or $1$ according as $\omega=\iota$ or $\varepsilon.$ 
For $T \in S_n({\bf Z}_p)_e,$ put $T^{(0)}=2^{-1}T,F_p^{(0)}(T,X)=F_p(T^{(0)},X),$ and $\widetilde F_p^{(0)}(T,X)=\widetilde F_p(T^{(0)},X).$  
  For $d_0 \in {\mathcal F}_{p}$ and  $\omega=\varepsilon^l$ with $l=0,1,$  we  define a  formal power series $P_{n,p}^{(0)}(d_0,\omega,X,t)$ in $t$ by 
$$P_{n,p}^{(0)}(d_0,\omega,X,t)=\kappa(d_0,n,l)^{-1} \hspace*{-2.5mm}\sum_{B \in {\mathcal L}_{n,p}^{(0)}(d_0)}  \hspace*{-2.5mm}{\widetilde F_p^{(0)}(B,X) \over \alpha_p(B)}\omega(B)t^{\nu_p(\det B)},$$
where  $$\kappa(d_0,n,l)=\kappa(d_0,n,l)_p =\{(-1)^{n(n+2)/8}((-1)^{n/2}2,d_0)_2\}^{l\delta_{2,p}}.$$
Let ${\mathcal F}$ denote the set of fundamental discriminants, and for $l=\pm 1,$ put 
$${\mathcal F}^{(l)}=\{ d_0 \in {\mathcal F} \  | \ ld_0 >0 \}.$$

\bigskip

\noindent
 {\bf Theorem 4.2.}  {\it 
 Let the notation and the assumption be as above.  Then for ${\rm Re}(s) \gg 0,$ we have 
 \begin{eqnarray*}
 \lefteqn
{L^*(s,I_n(h),\chi) =\kappa_{n}2^{ns-1-n/2}}
\\
& \times & \hspace*{-1mm} 
\{ \sum_{d_0 \in {\mathcal F}^{((-1)^{n/2})} }c_h(|d_0|)|d_0|^{n/4-k/2+1/4} \prod_p   P_{n,p}^{(0)}(d_0,\iota_p,\beta_p,p^{-s+k/2+n/4+1/4}\chi(p))  \\
&+ & (-1)^{n(n+2)/8} \\
& \times &\sum_{d_0 \in {\mathcal F}^{((-1)^{n/2})} }((-1)^{n/2} 2, d_0)_2c_h(|d_0|)|d_0|^{n/4-k/2+1/4} \prod_p P_{n,p}^{(0)}(d_0,\varepsilon_p,\beta_p,p^{-s+k/2+n/4+1/4}\chi(p)) \}. 
\end{eqnarray*}
} 

\bigskip
\begin{proof}
Let $T \in {{\mathcal L}_n}_{>0}.$ Then  the $T$-th Fourier coefficient $c_{I_n(h)}(T)$  of $I_n(h)$ is uniquely determined by the genus to which $T$ belongs, and, by definition, it can be expressed as
$$c_{I_n(h)}(T)=c_h(|{\textfrak d}_T|){\textfrak f}_T^{k-n/2-1/2}\prod_p \widetilde F^{(0)}(2T,\beta_p)$$
We also note that 
$$  {\textfrak f}_T^{k-n/2-1/2}=|{\textfrak d}_T|^{-(k/2-n/4-1/4)} (\det (2T))^{(k/2-n/4-1/4)} $$
and $e(2T)=e(T)$ for  $T \in {{\mathcal L}_n}_{>0}.$
Hence we have
$$\sum_{2T' \in {\mathcal G}(2T)} {c_{I_n(h)}(T') \over e(2T')}=\det (2T)^{k/2+n/4+1/4} |{\textfrak d}_T|^{-k/2+n/4+1/4}c_h(|{\textfrak d}_T|)\prod_p {\widetilde F_p^{(0)}(2T,\beta_p) \over \alpha_p(2T)}.$$
Thus, similarly to \cite{I-K1}, Theorem 3.3, (1), and \cite{I-K2}, Theorem 3.2,  we obtain
$$L^*(s,I_n(h),\chi) = \kappa_{n}2^{ns-1-n/2} \sum_{d_0 \in {\mathcal F}^{((-1)^{n/2})}} c_h(|d_0|)|d_0|^{n/4-k/2+1/4} $$
$$ \times \{\prod_p  P_{n,p}^{(0)}(d_0,\iota_p,\beta_p,p^{-s+k/2+n/4+1/4}\chi(p))$$
$$+  (-1)^{n(n+2)/8}((-1)^{n/2} 2, d_0)_2 \prod_p  P_{n,p}^{(0)}(d_0,\varepsilon_p,\beta_p,p^{-s+k/2+n/4+1/4}\chi(p)) \}.$$ 
This proves the assertion. 

\end{proof}

\bigskip

\noindent
{\bf Proposition 4.3.}     {\it 
Let $d_0 \in {\mathcal F}_{p}$ and $\xi_0=\widetilde \xi(d_0).$  
 Then 
  $$P_n^{(0)}(d_0,\iota,X,t)={(p^{-1}t)^{\nu_p(d_0)} \over \phi_{n/2-1}(p^{-2})(1-p^{-n/2}\xi_0)}$$
  $$\times  {(1+t^2p^{-n/2-3/2})(1+t^2p^{-n/2-5/2}\xi_0^2)-
  \xi_0 t^2p^{-n/2-2}(X+X^{-1}+p^{1/2-n/2} +p^{-1/2+n/2}) \over 
  (1-p^{-2}Xt^2)(1-p^{-2}X^{-1}t^2)\prod_{i=1}^{n/2} (1-t^2p^{-2i-1}X)(1-t^2p^{-2i-1}X^{-1})  },$$
  and
  $$P_n^{(0)}(d_0,\varepsilon,X,t)={1  \over \phi_{n/2-1}(p^{-2})(1-p^{-n/2}\xi_0)}{\xi_0^2 \over 
  \prod_{i=1}^{n/2} (1-t^2p^{-2i}X)(1-t^2p^{-2i}X^{-1})  }.$$     
  }
\begin{proof} Put $H_k=\mattwo{O}{1_k}{1_k}{O},$ and for $d \in {\bf Z}_p^*$  put
$$D=\{x \in M_{2k,n}({\bf Z}_p) \ | \ \det (H_k[x]) \in dp^i { {\bf Z}_p^*}^{\Box} \ {\rm with \ some} \ i \in {\bf Z}_{\ge 0} \}.$$
We then define $Z_{2k}(u,\varepsilon^l,d)$ as
$$Z_{2k}(u,\varepsilon^l,d)=\int_D \varepsilon^l(H_k[x]) |\det (H_k[x])|_p^{s-k} dx$$
with $u=p^{-s},$ where $|*|_p$ denotes the normalized valuation on ${\bf Q}_p,$ and $dx$ is the measure on $M_{2k,n}({\bf Q}_p)$ normalized so that the volume of $M_{2k,n}({\bf Z}_p)$ is $1.$ Moreover put
$$Z_{2k,e}(u,\varepsilon^l,d)={1  \over 2}(Z_{2k,n}(u,\varepsilon^l,d)+Z_{2k,n}(-u,\varepsilon^l,d)),$$
and
$$Z_{2k,o}(u,\varepsilon^l,d)={1  \over 2}(Z_{2k,n}(u,\varepsilon^l,d)-Z_{2k,n}(-u,\varepsilon^l,d)).$$
Then it is well known that
$$Z_{2k,x(d_0) }(u,\varepsilon^l,(-1)^{n/2}p^{-\nu_p(d_0)} d_0)=\phi_{n}(p^{-1})\sum_{T \in {\mathcal L}_{n,p}^{(0)}(d_0)} {b_p(2^{-\delta_{2,p}}T,p^{-k}) \over \alpha_p(T)}(p^k t)^{\nu_p(\det (T)}$$
for $d_0 \in {\mathcal F}_p,$ where $x(d_0)=e$ or $o$ according as $\nu_p(d_0)$ is even or odd. 
Recall that 
$$b_p(2^{-\delta_{2,p}}T,p^{-k})={(1-p^{-k})\prod_{i=1}^{n/2}(1-p^{-2k+2i}) \over 1-\xi(2^{-\delta_{2,p}}T)p^{-k+n/2}} F_p^{(0)}(T,p^{-k})$$
and 
$$F_p^{(0)}(T,p^{-k})= p^{(-k/2+(n+1)/4)(\nu_p(\det T)-\nu_p(d_0)) } \widetilde F_p^{(0)}(T,p^{-k+(n+1)/2}).$$
Hence we have
$$Z_{2k,x(d_0)}(u,\varepsilon^l,(-1)^{n/2}p^{-\nu_p(d_0)} d_0)=\phi_{n}(p^{-1}){(1-p^{-k})\prod_{i=1}^{n/2}(1-p^{-2k+2i}) \over 1-\xi(2^{-\delta_{2,p}}T)p^{-k+n/2} }$$
$$ \times  p^{(k/2-(n+1)/4)\nu_p(d_0)}P_n^{(0)}(d_0,\varepsilon^l,p^{-k+(n+1)/2},up^{k/2+(n+1)/4}).$$
Let $T(d_0,\omega,X,t)$ denote the right-hand side of the formula for $\omega=\varepsilon^l \ (l=0,1)$ in the proposition. Then, by [\cite{Sai2}, Theorem 3.4 (2)], we have 
$$Z_{2k,x(d_0)}(u,\varepsilon^l,(-1)^{n/2}p^{-\nu_p(d_0)} d_0)=\phi_{n}(p^{-1}){(1-p^{-k})\prod_{i=1}^{n/2}(1-p^{-2k+2i}) \over 1-\xi(T)p^{-k+n/2} }$$
$$ \times  p^{(k/2-(n+1)/4)\nu_p(d_0)}T(d_0,\varepsilon^l,p^{-k+(n+1)/2},up^{k/2+(n+1)/4}).$$
(Remark that there are misprints in \cite{Sai2}; the $(q^{-1})_n$ on page 197, lines 9 and 15 should be $(q^{-1})_r.$)
Hence we have 
$$P_n^{(0)}(d_0,\varepsilon^l,p^{-k+(n+1)/2},up^{k/2+(n+1)/4})=T(d_0,\varepsilon^l,p^{-k+(n+1)/2},up^{k/2+(n+1)/4})$$
for infinitely many positive integers $k$. Hence we have 
$$P_n^{(0)}(d_0,\varepsilon^l,X,t)=T(d_0,\varepsilon^l,X,t).$$
\end{proof}

{\bf Proof of Theorem 4.1.}  

Put $\Omega=\{\omega_p \},$ and let $d_0 \in {\mathcal F}^{((-1)^{n/2})}.$ Put
$$P(s,d_0,\Omega,\chi)=\prod_p  P_{n,p}^{(0)}(d_0,\varepsilon_p,\beta_p,p^{-s+k/2+n/4+1/4}\chi(p))
.$$
Then by Proposition 4.3,  we have
$$P(s,d_0,\{\iota_p \},\chi)$$
$$=|d_0|^{-s+k/2+n/4-3/4}\chi(d_0) \prod_{i=1}^{n/2-1} \zeta(2i) L(n/2,\chi_{d_0})  \prod_{i=1}^{n/2} L(2s+2i-n,S(h),\chi^2)$$
$$ \times L(2s-n+1,S(h),\chi^2)\prod_p  \{(1+p^{-2s+k-1}\chi(p)^2)(1+\chi_{d_0}(p)^2 p^{-2s+2k-2}\chi(p)^2)$$
$$-\chi_{d_0}(p) p^{-2s+k-3/2}\chi(p)^2\beta_p (1+p^{1/2-n/2}\beta_p^{-1})(1+p^{-1/2+n/2}\beta_p^{-1}) \}.$$
We note that $L(s, h)$ and $L(s,E_{n/2+1})$ can be expressed as
$$L(s,h)=L(2s,S(h))\sum_{d_0 \in {\mathcal F}^{((-1)^{n/2})}} c(|d_0|)|d_0|^{-s} \prod_p (1-\chi_{d_0}(p)p^{k-n/2-1-2s}),$$
and  
$$L(s,E_{n/2+1})=\zeta(2s)\zeta(2s-n+1)$$
$$ \times \sum_{d_0 \in {\mathcal F}^{((-1)^{n/2})}} L(1-n/2,\chi_{d_0})|d_0|^{-s} \prod_p (1-\chi_{d_0}(p)p^{n/2-1-2s}),$$
and therefore,  we easily see that $L(s, h,E_{n/2+1/2},\chi)$ can be expressed as 
$$R(s,h,E_{n/2+1/2},\chi)= L(2s,S(h),\chi^2)L(2s-n+1,S(h),\chi^2)$$
$$\times \sum_{d_0 \in {\mathcal F}^{((-1)^{n/2})}}|d_0|^{-s} c(|d_0|)\chi(d_0)L(1-n/2, \chi_{d_0})$$
$$ \times \prod_p \{ (1+p^{-2s+k-1}\chi(p)^2)(1+\chi_{d_0}(p)^2 p^{-2s+k-2}\chi(p)^2)$$
$$ -\chi_{d_0}(p) p^{-2s+k-3/2}\chi(p)^2 \beta_p (1+p^{1/2-n/2}\beta_p^{-1})(1+p^{-1/2+n/2}\beta_p^{-1})\}$$
(cf. \cite {Sh1}, Lemma 1.)
Thus, by remarking the functional equation 
$$L(1-n/2,\chi_{d_0})=2^{1-n/2}\pi^{-n/2}{\rm \Gamma}(n/2)|d_0|^{(n-1)/2}L(n/2,\chi_{d_0}),$$
 we have
$$\sum_{d_0 \in {\mathcal F}^{((-1)^{n/2})}}c_h(|d_0|)|d_0|^{-s+k/2+n/4+1/4}  P(s,d_0,\{\iota_p \},\chi)$$
$$=\prod_{i=1}^{n/2-1} \zeta(2i){2^{n/2-1}\pi^{n/2}   \over {\rm \Gamma}(n/2)}  R(s,h,E_{n/2+1/2};\chi)  \prod_{i=1}^{n/2-1} L(2s-2i+n,S(h),\chi^2).$$
On the other hand, if $d_0 \not=1,$ by Proposition 4.3,
 we have 
$$P(s,d_0,\{\varepsilon_p \},\chi)=0.$$
Thus if $n \equiv 2 \ {\rm mod} \ 4,$ for any $d_0 \in {\mathcal F}^{((-1)^{n/2})},$  $$P(s,d_0,\{\varepsilon_p\},\chi)=0.$$
If $n \equiv 0 \ {\rm mod} \ 4,$ by Proposition 4.3, we have
$$ P(s,1,\{\varepsilon_p\},\chi)=\zeta(n/2)\prod_{i=1}^{n/2-1}\zeta(2i)\prod_{i=1}^{n/2} L(2s-2i+1,S(h),\chi^2).$$
Thus the assertion follows from Theorem 4.2. \hfill $\Box$

\section{Relation between twisted Kocher-Maa{\ss} series of the first and second kinds}
Let $N$ be a positive integer. Let $g$ be a periodic function on ${\bf Z}$ with a period $N$ and $\phi$ a polynomial in $t_1,\cdots ,t_r.$ Then for an element $u=(a_1 \ {\rm mod} \ N, \cdots ,a_r \ {\rm mod} \ N) \in ({\bf Z}/N{\bf Z})^r,$ the value $g(\phi(a_1, \cdots ,a_r))$ does not depend on the choice of the representative of $u.$ Therefore we denote this value by $g(\phi(u)).$ In particular we sometimes regard a Dirichlet character mod $N$ as a function on ${\bf Z}/N{\bf Z}.$ 

For a Dirichlet character $\chi$ mod $N$ and $A \in {\mathcal L}_{m>0},$ put
$$h(A,\chi)=
\sum_{ U \in SL_{m,N}({\bf Z}) \backslash SL_m({\bf Z})} \chi({\rm tr}(A[U])).$$
 As was shown in [\cite{K-M}, Proposition 3.1], the twisted Koecher-Maa{\ss} series of the first kind of a Siegel modular form can be expressed in terms of $h(A,\chi)$ as stated later. Therefore we shall compute $h(A,\chi)$ in the case where $A$ is an element of ${\mathcal L}_{m>0}.$ 
 For $A=(a_{ij})_{m \times m} \in S_m({\bf Z}/N{\bf Z})$ and $c \in {\bf Z}/N{\bf Z},$ put
 $${\mathcal R}_{N}(A,c)=\{X=(x_{ij})_{m \times m} \in M_m({\bf Z}/N{\bf Z}) \ | \ \sum_{i=1}^m \sum_{\alpha,\beta=1}^m a_{\alpha,\beta}x_{i \alpha} x_{i \beta}- c = 0 $$ 
$$ \ {\rm and} \ \det X -1 = 0 \}.$$
Then we have 
$$h(A,\chi)=\sum_{c \in {\bf Z}/N{\bf Z}}\chi(c) \#(R_N(A,c)).$$ 
From now on let $p$ be an odd prime number and $F_p$ be the field with $p$ elements. For $S \in S_m(F_p)$ and $T \in S_r(F_p)$ put 
$${\rm A}(S,T)=\{Y \in M_{r,m}(F_p) \ | \ YS\ {}^t \! Y=T\}.$$ For an element $S \in S_m(F_p)$ with $m$ even
 put $\chi(S)=\left({(-1)^{m/2}  \det S \over p} \right).$ 

\bigskip

\noindent
{\bf Lemma 5.1.} {\it
 Let $S \in S_m(F_p)^{\times}.$ \\
{\rm (1)} Let $T \in S_r(F_p)^{\times}$ with $m \ge r.$ \\
{\rm (1.1)} Let $r$ be even. Then
$$\#{\rm A}(S,T)=p^{rm-r(r+1)/2}(1-\chi(S)p^{-m/2})(1+\chi((-S) \bot T)p^{(r-m)/2})\prod_{m-r+1 \le e \le m-1 \atop e \ {\rm even} } (1-p^{-e})$$
or 
$$\#{\rm A}(S,T)=p^{rm-r(r+1)/2}\prod_{m-r+1 \le e \le m-1 \atop e \ {\rm even} } (1-p^{-e})$$
according as $m$ is even or odd. \\
{\rm (1.2)} Let $r$ be odd. Then
$$\#{\rm A}(S,T)=p^{rm-r(r+1)/2}(1-\chi(S)p^{-m/2})\prod_{m-r+1 \le e \le m-1 \atop e \ {\rm even} } (1-p^{-e})$$
or
$$\#{\rm A}(S,T)=p^{rm-r(r+1)/2}(1+\chi((-S) \bot T)p^{(r-m)/2})\prod_{m-r+1 \le e \le m-1 \atop e \ {\rm even} } (1-p^{-e})$$
according as $m$ is even or odd. In particular, for $c \in F_p^{\times},$ we have
$$\#{\rm A}(S,c)=p^{m/2-1}(p^{m/2}-\left({(-1)^{m/2} \det S \over p} \right))$$
or
$$\#{\rm A}(S,c)=p^{(m-1)/2}(p^{(m-1)/2}+\left({(-1)^{(m+1)/2} c \det S \over p} \right))$$
according as $m$ is even or odd. \\
{\rm (2)} We have
$$\#{\rm A}(S,0)=p^{m/2-1}(p^{m/2}-\left({(-1)^{m/2} \det S \over p} \right))+p^{m/2}\left({(-1)^{m/2} \det S \over p} \right)$$
or 
$$\#{\rm A}(S,0)=p^{m-1}$$
according as $m$ is even or odd.
 }
\begin{proof}
The assertions (1) and (2) follow from [\cite{Ki1}, Theorem 1.3.2], and  [\cite{Ki1}, Lemma 1.3.1], respectively.  
\end{proof}

\bigskip

\noindent
{\bf Proposition 5.2.} {\it
 Let $A=a_1 \bot \cdots \bot a_m$ with $a_i \in F_p.$ 
For $c \in F_p^{\times}$ put 
$${\mathcal M}_p(A,c)=\{Z=(z_{ij}) \in S_{m}(F_p) \ | \ \det (Z) =1 \ {\rm and} \ c-\sum_{i=1}^m a_iz_{ii}=0 \},$$
and
$$\gamma_{m,p}=p^{m^2-m(m+1)/2}(1-p^{-m/2})\prod_{e=1}^{(m-2)/2} (1-p^{-2e})$$
or 
$$\gamma_{m,p}=p^{m^2-m(m+1)/2}\prod_{e=1}^{(m-1)/2} (1-p^{-2e})$$
according as $m$ is even or odd. Then we have
$$\#{\mathcal R}_p(A,c)=\gamma_{m,p}\#{\mathcal M}_p(A,c) .$$
}

\begin{proof} Let $\Phi:GL_m(F_p) \longrightarrow S_m(F_p) \cap GL_m(F_p)$ be the mapping defined by $\Phi(X)=X{}^t \! X.$
Then by Lemma 5.1, we have $\#\Phi^{-1}(Z)=2\gamma_{m,p}$ for any $Z \in S_m(F_p) \cap SL_m(F_p).$
We note that $\det X =\pm 1$ for any $X \in \Phi^{-1}(Z).$ Hence we have $\#(\Phi^{-1}(Z) \cap SL_m(F_p))=\gamma_{m,p}.$
Moreover we have 
$${\rm tr}({}^t \! X AX)={\rm tr}(AX{}^t \! X),$$
and hence $X \in {\mathcal R}_p(A,c)$ if and only if $\Phi(X) \in {\mathcal M}_p(A,c).$ This proves the assertion. 
\end{proof}

We rewrite ${\mathcal M}_p(A,c)$ in more concise form. 
For a positive integer $N$ we denote by ${\mathcal D}_N$ the set of Dirichlet charcters mod $N,$ and for a positive integer $m$ we denote by ${\mathcal D}_{N,m}$ the subset of ${\mathcal D}_N$ consisiting of Dirichlet characters whose $m$-th power is the trivial character.  
We note that ${\mathcal D}_{p,m}={\mathcal D}_{p,l}$ with $l={\rm GCD}(m,p-1)$ if $p$ is an odd  prime number. 
 We denote by $\displaystyle \left({* \over N}\right)$ the Jacobi symbol for a positive odd integer $N$. For two Dirichlet characters $\chi$ and $\eta$ mod $N$ we define $J_m(\chi,\eta)$ and $I_m(\chi,\eta)$ 
$$J_m(\chi,\eta)=\sum_{Z \in S_m({\bf Z}/N{\bf Z})} \chi(\det Z)\eta(1-{\rm tr} (Z))$$
and
$$I_m(\chi,\eta)=\sum_{Z \in S_m({\bf Z}/N{\bf Z})} \chi(\det Z)\eta({\rm tr}(Z)).$$
By definition, $J_m(\chi,\eta)$ is an algebraic number. We note that $J_1(\chi,\eta)$ is the Jacobi sum
 $J(\chi,\eta)$ associated with $\chi$ and $\eta.$  We also define $J_m(\chi)$ as 
$J_m(\chi)= J_{m}({\chi}\left({ * \over N}\right)^{m-1} , {\chi}).$

\bigskip

\noindent
{\bf Lemma 5.3.} {\it 
Let $\eta$ be a primitive character mod $p.$ Let $c \in F_p$ and  $S \in S_l(F_p)$ of rank $r.$ Let $S \sim S_0 \bot O_{l-r}$ with $\det S_0 \not=0.$   Put
$$I_{\eta,S,c}=\sum_{{\bf w} \in F_p^l} \eta(S[^{t}{\bf w}]+c).$$
Assume that $r$ is odd, and that $\eta^2 \not=1.$  Then 
$$I_{\eta,S,c}=p^{l-(r+1)/2}J(\eta,({* \over p}))\left({(-1)^{(r+1)/2} \det S_0 \over p}\right)\eta(c) \left({c \over p}\right).$$
Assume that $r$ is even, and that $\eta \not=1.$  Then
$$I_{\eta,S,c}=p^{l-r/2} \left({(-1)^{r/2} \det S_0 \over p}\right)\eta(c).$$
Here we make the convention that $\left({(-1)^{r/2} \det S_0 \over p}\right)=1$ if $r=0.$
}

\begin{proof}
We have
$$I_{\eta,S,c}=p^{l-r} I_{\eta,S_0,c}.$$
Hence we may assume that $r=l.$ Then
$$I_{\eta,S,c}=\sum_{u \in F_p} \eta(u)\#A(S,u-c).$$
Let $l$ be odd. Then by Lemma 5.1, 
$$\#A(S,u-c)=p^{(l-1)/2}(p^{(l-1)/2}+\left({(-1)^{(l-1)/2} (u-c)\det S\over p}\right)).$$
Hence we have 
$$I_{\eta,S,c}= p^{(l-1)/2}\left({(-1)^{(l+1)/2} \det S \over p}\right)\sum_{u \in F_p} \eta(u) \left({u-c \over p}\right).$$
Since $\eta^2$ is nontrivial, we have $I_{\eta,S,c}=0$ if $c=0.$ If $c \not=0,$ then
$$\sum_{u \in F_p} \eta(u) \left({u-c \over p}\right)=\left({-c \over p}\right)\sum_{u \in F_p} \eta(u) \left({1-c^{-1}u \over p}\right)$$
$$=\eta(c)\left({-c \over p}\right)\sum_{u \in F_p} \eta(v) \left({1-v \over p}\right) =\eta(c) \left({-c \over p}\right)J(\eta, \left({* \over p}\right)).$$
Let $l$ be even. Then
$$\#A(S,u-c)=(p^{l/2}-\left({(-1)^{l/2} \det S \over p}\right))p^{l/2-1}+ p^{l/2}\left({(-1)^{l/2} \det S \over p}\right)a_0 ,$$
where $a_0=1$ or $0$ according as $u=c$ or not. Hence 
$$I_{\eta,S,c}= p^{l/2}\left({(-1)^{l/2} \det S \over p}\right)\eta(c).$$
\end{proof}
{\bf Corollary.} {\it
Let $d \in F_p^{\times}.$ Then we have
$$I_{\eta,S,cd}= \eta(d) \left({d \over p}\right)^r  I_{\eta,S,c}.$$
}

\bigskip

\noindent
{\bf Proposition 5.4.} {\it
Let $\eta$ be a primitive character mod $p.$  For  $Z_1 \in S_{l-1}(F_p)$ and $z_{ll} \in F_p,$ put 
$$I(Z_1,z_{ll})=\sum_{w \in M_{l-1,1}(F_p)}\eta(\det \left(\begin{smallmatrix} Z_1 & w \\ {}^tw & z_{ll} \end{smallmatrix}\right) ).$$
{\rm (1)} Assume that $l$ is even, and that $\eta^2 \not=1.$ Then
$$I(Z_1,z_{ll})=p^{(l-2)/2}J(\eta,({* \over p}))\left({(-1)^{l/2}  \det Z_1 \over p}\right)\eta(\det Z_1 z_{ll}) \left({ z_{ll} \over p}\right).$$
{\rm (2)} Assume that $l$ is odd, and that $\eta^2 \not=1.$ Then
$$I(Z_1,z_{ll})=p^{(l-1)/2}\left({(-1)^{(l-1)/2} \det Z_1 \over p}\right)\eta(\det Z_1 z_{ll}).$$
}
\begin{proof}
We note that
$$\det \mattwo {Z_1}{w}{{}^tw}{ z_{ll}} =-{\rm Adj}(Z_1)[w]+ \det Z_1 z_{ll},$$
where ${\rm Adj}(Z_1)$ is the $(l-1) \times (l-1)$ matrix whose $(i,j)$-th component is the 
 $(j,i)$-th cofactor of $Z_1.$ We also note that 
 $\det (-{\rm Adj}(Z_1))= (-1)^{l-1}(\det Z_1)^{l-2}.$ Thus the assertion follows directly from Lemma 5.3 if $\det Z_1 \not=0.$
 If $\det Z_1=0,$ then ${\rm rank}_{F_p} {\rm Adj}(Z_1) \le 1,$  and the assertion follows also from Lemma 5.3. 
 \end{proof}
 
 \bigskip
 
 Let $\chi$ be a Dirichlet character of odd conductor and  $A \in {\mathcal L}_{m>0}.$ Then $2^{2[m/2]}\det A$ belongs to ${\bf Z},$ and we define $\chi(\det A)$ as $\overline {\chi(2^{2[m/2]})} \chi(2^{2[m/2]}\det A).$  
 
\bigskip

 \noindent
{\bf Theorem 5.5.} {\it 
Let $\chi$ be a primitive character mod $p.$ Let $l={\rm GCD}(m,p-1),$ and $u_0$ be a primitive $l$-th root of unity mod $p.$ Let  $A \in {\mathcal L}_{m>0}.$ \\
{\rm (1)} If $\chi(u_0) \not=1,$ then we have $h(A,\chi)=0.$ \\
{\rm (2)} Assume that $\chi(u_0)=1.$ Fix  a  character $\widetilde \chi$  such that $\widetilde \chi^m=\chi.$\\
{\rm (2.1)} Let $m$ be even. Then 
$$h(A,\chi)=\gamma_{m,p} A_{m,p}\sum_{\eta \in {\mathcal D}_{p,m}}  (\widetilde \chi \eta)(\det A) J(\overline{\widetilde \chi \eta }, \left({* \over p}\right))J_{m-1}(\overline {\widetilde \chi \eta }),$$
where
$A_{m,p}=p^{(m-2)/2}(-1)^{m(p-1)/4 }.$ \\
{\rm (2.2)} Let $m$ be odd and assume that $\chi^2 \not=1.$ Then 
$$h(A,\chi)=\gamma_{m,p} A_{m,p}\sum_{ \eta \in {\mathcal D}_{p,m}}  (\widetilde \chi \eta)(\det A) J_{m-1}(\overline {\widetilde \chi \eta}),$$
where $A_{m,p}=p^{(m-1)/2}(-1)^{(m-1)(p-1)/4}.$ 
 }

\noindent
{\bf Remark.} The above formulation is based on the referee's suggestion. In the original version, we formulated Theorems 5.5 and 5.6 in terms of ``modified  power residue symbols".
\begin{proof}
We may regard $A$ as an element of $S_{m}(F_p).$ 
If $A=O_m$ then we have $h(A,\chi)=0.$ Hence we  assume that $A \not=O_m.$ Then we may assume that
 $A=a_1 \bot \cdots \bot a_{m-1} \bot  d$ with $d \not=0.$ 
Put 
$$\widetilde {\mathcal M}_p(A,c)$$
$$=\{ (Z_1,w) \in S_{m-1}(F_p)  \times M_{m-1,1}(F_p) \ | \ \det \mattwo{Z_1}{w}{{}^tw}{d^{-1}(1-\sum_{i=1}^{m-1} a_iz_{ii})} c^m=1  \}.$$
Write $Z \in S_m(F_p)$ as $Z=\smallmattwo{Z_1}{w}{{}^tw}{z_m}$ with $Z_1 \in S_{m-1}(F_p), w \in M_{m-1,1}(F_p),z \in F_p.$ 
Then the mapping $S_m(F_p) \ni Z \mapsto (c^{-1}Z_1,c^{-1}w)  \in S_{m-1}(F_p)  \times M_{m-1,1}(F_p) $ induces a bijection from
 ${\mathcal M}_p(A,c)$ to $\widetilde {\mathcal M}_p(A,c),$ and hence $\#\widetilde {\mathcal M}_p(A,c)=\#{\mathcal M}_p(A,c).$
Put
$$K(A)=\sum_c \#\widetilde {\mathcal M}_p(A,c)\chi(c).$$
Assume that $\chi(u_0) \not=1.$ Then 
we have
$$K(A)=\sum_{c \in F_p} \chi(cu_0) \# \widetilde  {\mathcal M}_p(A,cu_0).$$
We note that $\widetilde  {\mathcal M}_p(A,cu_0)=\widetilde  {\mathcal M}_p(A,c).$ Hence we have
$$K(A)=\chi(u_0)K(A).$$
Hence we have $K(A)=0.$ 

Assume that $\chi(u_0) =1.$ Then we can take  a Dirichlet character $\widetilde \chi$ such that $\widetilde \chi^m=\chi.$ 
First assume that $\det A=0.$ Then we may assume that we have
$A=A_0 \bot 0$ with $A_0 \in S_{m-1}(F_p).$ 
Let $P_{m-1,m}$ be the set of $(m-1) \times m$ matrices with entries in $F_p$ of rank $m-1.$ Then for each $X_1 \in P_{m-1,m}$ there exist exactly 
$p^{m-1}$ elements $X_2 \in M_{1,m}(F_p)$ such that $\vecttwo{X_1}{X_2} \in SL_m(F_p).$
Hence we have
$$h(A,\chi)=p^{m-1}\sum_{X_1 \in P_{m-1,m}} \chi(A_0[X_1]).$$
Let $m$ be even. Then we can take an element $\alpha \in F_p^{\times}$ such that $\chi(\alpha) \not=1.$ 
Moreover we can take $U_0 \in GL_{m}(F_p)$ such that 
${}^tU_0U_0=\alpha 1_m$ in view of (1.1) of Lemma 5.1. Hence 
$$h(A,\chi)=p^{m-1}\sum_{X_1 \in P_{m-1,m}} \chi(A_0[X_1U_0])=\chi(\alpha)h(A,\chi).$$
Hence we have $h(A,\chi)=0.$ Let $m$ be odd and assume that $\chi^2 \not=1.$ Then we can take an element $\alpha \in (F_p^{\times})^{\Box}$ such that $\chi(\alpha) \not=1.$ 
Moreover we can take $U_0 \in GL_{m}(F_p)$ such that ${}^tU_0U_0=\alpha 1_m$ in view of (1.2) of Lemma 5.1. 
Thus by the same argument as above we have $h(A,\chi)=0.$
This proves the assertion. 
Next assume that $\det A \not =0.$ We may assume that 
$$A=1_{m-1}  \bot d$$
with $d =\det A.$
Then we have
$$K(A)=\sum_c \#\widetilde {\mathcal M}_p(A,c)\widetilde \chi(c^m).$$
Hence we have
$$K(A)=\sum_{(Z_1,w) } \overline{ \widetilde \chi(\det {\smallmattwo{Z_1}{w}{{}^tw}{d^{-1}(1-{\rm tr} (Z_1))}}) },$$
where $(Z_1,w)$ runs over elements of $S_{m-1}(F_p) \times  M_{m-1,1}(F_p)$ such that 
$$(*) \hskip 3cm \det \mattwo{Z_1}{w}{{}^tw} {d^{-1}(1-{\rm tr} (Z_1))} =u^m$$
 with some $u \in F_p^{\times},$ and for such a matrix $\mattwo{Z_1}{w}{{}^tw} {d^{-1}(1-{\rm tr} (Z_1))},$ there exist exactly
 $l$ elements $u$ of $F_p$ satisfying (*).
We have
$$\sum_{\eta \in {\mathcal D}_{p,m}} \eta(v) =l \ {\rm or} \ 0$$
according as $v =u^m$ with some $u \in F_p^{\times}$ or not. 
Hence we have 
$$K(A)=\sum_{\eta \in {\mathcal D}_{p,m}} \sum_{(Z_1,w) } \overline{ (\widetilde \chi \eta)(\det {\smallmattwo{Z_1}{w}{{}^tw}{ d^{-1}(1-{\rm tr} (Z_1))}})}.$$
For each $\eta \in {\mathcal D}_{p,m}$ put
$$K(A)_{\eta} = \sum_{(Z_1,w) }  \overline{ (\widetilde \chi  \eta)(\det {\smallmattwo{Z_1}{w}{{}^tw}{d^{-1}(1-{\rm tr} (Z_1))}})
}$$
We note that $(\widetilde \chi \eta)^2 \not=1$ for any $\eta.$  Hence by Proposition 5.4 we have
$$K(A)_{\eta} =A_{m,p} J_{\eta}\sum_{Z_1 \in S_{m-1}(F_p)  } (\widetilde \chi \eta)^*(\det A)\overline {(\widetilde \chi \eta)^*(\det Z_1)} \overline{(\widetilde \chi \eta)^*(1-{\rm tr}(Z_1)) } \left({ \det Z_1 \over p}\right)^m, 
$$
where $J_{\eta}= J(\overline{\widetilde \chi \eta }, \left({* \over p}\right))$ or $1$ according as $m$ is even or odd, and $(\widetilde \chi \eta)^*=\widetilde \chi \eta \left({* \over p}\right)^{m-1}.$ This proves the assertion if $m$ is odd. Assume that $m$ is even. Then it is easily seen that the set $\{\widetilde \chi \eta \left({* \over p}\right)\}_{\eta \in {\mathcal D}_{p,m}}$ of Dirichlet characters coincides with $\{ \widetilde \chi \eta \}_{\eta \in {\mathcal D}_{p,m}}.$ Moreover $(\widetilde \chi \eta)^2 \not=1$ for any $\eta.$ 
This proves the assertion. 
\end{proof}

\bigskip
Let $\chi$ be a Dirichlet character mod $N.$ Fix a prime factor $p$ of $N.$ For an integer $n$ prime to $p,$ take an integer $m$ such that 
$$m \equiv \listtwo{n \ {\rm mod} \  \ p^e}{1 \ {\rm mod} \ N/p^e}.$$
We then put  
$$\chi^{(p)}(n)= \listtwotwo{\chi(m) \ }{{\ {\rm if} \ (n,p)=1}}{0 \ }{{\ {\rm if} \ (n,p)\not=1}}.$$
Then it is independent of the choice of $m,$ and $\chi^{(p)}$ is a character mod $p^e,$ and we have $\chi=\prod_{p|N} \chi^{(p)}.$ 

\noindent
{\bf Theorem 5.6.} {\it
Let $N=p_1 \cdots p_r$ with $p_1,\cdots,p_r$ distinct odd prime numbers. Put $l_i={\rm G.C.D}(p_i-1,m).$ Let $\chi$ be a primitive Dirichlet character mod $N.$ Let $u_{0,i}$ be a  primitive $l_i$-th root of unity mod $p_i.$ Let $A \in {\mathcal L}_{m>0}.$ \\
{\rm (1)} If $\chi^{(p_i)}(u_{0,i}) \not=1$ for some $i.$ Then we have $h(A,\chi)=0.$ \\
{\rm (2)} Assume that $\chi^{(p_i)}(u_{0,i})=1$ for all $i.$ Fix a character $\tilde \chi$ such that  $\tilde \chi^m=\chi.$ \\
{\rm (2.1)} Let $m$ be even. Then we have 
$$h(A,\chi)= \prod_{i=1}^r (-1)^{m(p_i-1)/4 } p_i^{(m-2)/2}\gamma_{m,p_i} $$$$ \times \sum_{ \eta \in {\mathcal D}_{N,m}} ( \widetilde \chi \eta)(\det  A)J(\widetilde {\chi} \eta,  \left({* \over N}\right))  J_{m-1}(\overline {\widetilde  \chi \eta}).$$\\
{\rm (2.2)} Let $m$ be odd, and assume that $\chi^2$ is primitive. Then we have 
$$h(A,\chi)= \prod_{i=1}^r (-1)^{(m-1)(p_i-1)/4 }p_i^{(m-1)/2} \gamma_{m,p_i}  $$
$$ \times \sum_{ \eta \in {\mathcal D}_{N,m}} (\widetilde \chi \eta)(\det  A)J_{m-1}(\overline {\widetilde  \chi \eta}) .$$
}

\begin{proof}
 We note that the mapping ${\mathcal M}_N \ni \chi \mapsto (\chi^{(p_1)},\cdots,\chi^{(p_r)}) \in {\mathcal M}_{p_1} \times \cdots \times {\mathcal M}_{p_r}$ 
 induces a bijection from ${\mathcal M}_{N,m}$ to ${\mathcal M}_{p_1,m} \times \cdots \times {\mathcal M}_{p_r,m}.$  We also  note that $J_m(\eta_1,\eta_2)=\prod_{i=1}^r J_m(\eta_1^{(p_i)},\eta_2^{(p_i)})$
  for primitive characters $\eta_1$ and $\eta_2$ mod $N.$ Moreover $\eta_j^2$ is primitive if and only if  ${\eta_j^{(p_i)}}^2 \not=1$ for any $1 \le i \le r.$ Thus the assertion follows from  Theorem 5.5 and [\cite{K-M}, Lemma 3.2].
 \end{proof}

Now we give explicit formulas for $J_m(\chi,\eta)$ and $I_m(\chi,\eta).$ 

\bigskip

\noindent
{\bf Proposition 5.7.}{\it
Let $\chi$ and $\eta$ be primitive characters mod $p.$  Assume that $\chi^2 \not=1.$ Put $c_m(\chi,\eta)=1$ or $0$ according as $\chi^{m}\eta=1$ or not. \\
{\rm (1)} Assume that $m$ is odd. Then
$$I_m(\chi,\eta)=c_m(\chi,\eta)\left({-1 \over p}\right)^{(m-1)/2}p^{(m-1)/2}(p-1)J_{m-1}(\chi\left({ *\over p}\right),\eta).$$
{\rm (2)} Assume that $m$ is even. Then
$$I_m(\chi,\eta)=c_m(\chi,\eta)\left({-1 \over p}\right)^{m/2}p^{(m-2)/2}(p-1)\chi(-1)J(\chi, \left({ *\over p}\right))J_{m-1}(\chi\left({ *\over p}\right),\eta).$$
}

\begin{proof}
By Proposition 5.4, we have
$$I_m(\chi,\eta)=I_m'   \times 
\left\{\begin{array}{ll}
p^{(m-1)/2} \left({(-1)^{(m-1)/2}  \over p}\right) & \ {\rm if} \ m \ {\rm  is \ odd} \\
p^{(m-2)/2} \left({(-1)^{(m-2)/2}  \over p}\right)J(\chi, \left({ *\over p}\right)) & \ {\rm if} \ m \ {\rm is \ even,}
\end{array}
\right.
$$ 
where $$I'_m=\sum_{z_{mm} \in F_p  \atop Z_1 \in S_{m-1}(F_p) ^{\times}}\chi(z_{mm}) \chi (\det Z_1)\left({\det Z_1 \over p}\right) \eta(z_{mm}+{\rm tr}(Z_1))\left({ z_{mm} \over p}\right)^{m-1}.$$
Then we have
$$I_m'=\sum_{z_{mm} \in F_p^{\times}  \atop Z_1 \in S_{m-1}(F_p) ^{\times}}\chi(z_{mm})\eta(z_{mm}) \chi (\det Z_1)\left({\det Z_1 \over p}\right) \eta(1+z_{mm}^{-1}{\rm tr}(Z_1))\left({ z_{mm} \over p}\right)^{m-1} .$$
Put $Y_1=-z_{mm}^{-1}Z_1.$ Then $\det Y_1=(-1)^{m-1}z_{mm}^{1-m} \det Z_1.$ 
Hence we have 
$$I_m'=\chi((-1)^{m-1})\left({(-1)^{m-1}  \over p}\right)$$
$$ \times \sum_{z_{mm} \in F_p^{\times} }\chi(z_{mm})^m \eta(z_{mm})\sum_{Y_1 \in S_{m-1}(F_p) ^{\times}} \chi (\det Y_1)\left({\det Y_1 \over p}\right) \eta(1-{\rm tr}(Y_1)).$$
We have 
$$\sum_{z_{mm} \in F_p^{\times} }\chi(z_{mm})^m \eta(z_{mm})=p-1 \ {\rm or} \ 0$$
according as $\chi^m \eta$ is trivial or not. This proves the assertion. 
\end{proof}

\bigskip

\noindent
{\bf Proposition 5.8.} {\it 
Let $\chi$ and $\eta$ be as in Proposition 5.7.  \\
{\rm (1)} Assume that $m$ is odd. Then
$$J_m(\chi,\eta)=\left({(-1)^{(m-1)/2} \over p}\right)p^{(m-1)/2}$$
$$ \times \{J(\chi,\chi^{m-1}\eta)J_{m-1}(\chi\left({ *\over p}\right),\eta)+\eta(-1)I_{m-1}(\chi\left({ *\over p}\right),\eta)\}.$$
{\rm (2)} Assume that $m$ is even. Then
$$J_m(\chi,\eta)=\left({-1 \over p}\right)^{m/2}p^{(m-2)/2}J(\chi, \left({ *\over p}\right))$$
$$ \times \{J(\chi \left({ *\over p}\right),\chi^{m-1}\left({ *\over p}\right)\eta)J_{m-1}(\chi \left({ *\over p}\right) ,\eta)+\eta(-1)I_{m-1}(\chi \left({ *\over p}\right) ,\eta)\}.$$
}
\begin{proof}
By Proposition 5.4, we have
$$J_m(\chi,\eta)=(J_m'+J_m'') \times 
\left\{ \begin{array}{ll}
p^{(m-1)/2}\left({(-1)^{(m-1)/2} \over p}\right) & \ {\rm if} \ m \ {\rm is \ odd} \\
p^{(m-2)/2}\left({(-1)^{(m-2)/2} \over p}\right) J(\chi, \left({ *\over p}\right)) & \  {\rm if} \ m \ {\rm is \ even,}
\end{array}
\right.
$$
where 
$$J_m'=\sum_{z_{mm} \in F_p, z_{mm} \not=1  \atop Z_1 \in S_{m-1}(F_p) ^{\times}}\left({ \det Z_1 \over p}\right)\left({z_{mm} \over p}\right)^{m-1}\chi(z_{mm}) \chi (\det Z_1) \eta(1-z_{mm}-{\rm tr}(Z_1)),$$
and 
 $$J_m''=\sum_{ Z_1 \in S_{m-1}(F_p) ^{\times}} \left({ \det Z_1 \over p}\right) \chi (\det Z_1) \eta(-{\rm tr}(Z_1)).$$
Then we have  $J_m''=\eta(-1)I_{m-1}(\chi\left({ *\over p}\right)  ,\eta).$ Moreover 
$$J_m'=\sum_{z_{mm} \in F_p, z_{mm} \not=1  \atop Z_1 \in S_{m-1}(F_p) ^{\times}}\chi(z_{mm}) \left({ \det Z_1 \over p}\right)\left({ z_{mm} \over p}\right)^{m-1}\chi (\det Z_1) $$
$$ \times \eta(1-z_{mm}) \eta(1-(1-z_{mm})^{-1} {\rm tr}(Z_1)).$$
Put $Y_1=(1-z_{mm})^{-1} Z_1.$ Then $\det Y_1=(1-z_{mm})^{1-m}\det Z_1.$ Hence we have
$$J_m'=\sum_{z_{mm} \in F_p } \chi(z_{mm}) \left({ z_{mm} \over p}\right)^{m-1}  \left({1- z_{mm} \over p}\right)^{m-1}   \chi(1-z_{mm})^{m-1} \eta(1-z_{mm}) $$
$$ \times \sum_{Y_1 \in S_{m-1}(F_p) ^{\times}}  \left({ \det Y_1 \over p}\right)\chi (\det Y_1)  \eta(1- {\rm tr}(Y_1)).$$
This proves the assertion. 
\end{proof}

\bigskip

\noindent
{\bf Theorem 5.9.} {\it
Let $\chi$ be a primitive character mod $p.$ \\
 {\rm (1)} Let $m$ be odd, and assume that $\chi^2 \not=1.$  \\
{\rm (1.1)} Assume that $\chi^m \not=1.$ Then 
$$J_m(\chi \left({ *\over p}\right)^i,\chi)=\left({-1 \over p}\right)^{(m-1)/2}p^{(m-1)/2}J(\chi\left({ *\over p}\right)^i ,\chi^{m})J_{m-1}(\chi\left({ *\over p}\right)^{i+1},\chi).$$
{\rm (1.2)} Assume that $\chi^m =1.$ Then 
$$J_m(\chi \left({ *\over p}\right)^i,\chi)=p^{m-1}\left({ -1 \over p}\right)^{i+1}J(\chi \left({ *\over p}\right)^{i+1}, \left({ * \over p}\right))J_{m-2}(\chi \left({ *\over p}\right)^i,\chi).$$
{\rm (2)} Let $m$ be even. \\
{\rm (2.1)} Assume that $\chi^m \left({ *\over p}\right)^{i +1}\not=1.$ Then 
$$J_m(\chi \left({ *\over p}\right)^i,\chi)=\left({-1 \over p}\right)^{(m-2)/2}J(\chi \left({ *\over p}\right)^{i}, \left({ *\over p}\right)) J(\chi\left({ *\over p}\right)^{i+1} ,\chi^{m}\left({ *\over p}\right)^{i+1})J_{m-1}(\chi\left({ *\over p}\right)^{i+1},\chi).$$
{\rm (2.2)} Assume that $\chi^m \left({ *\over p}\right)^{i +1}=1.$ Then 
$$J_m(\chi \left({ *\over p}\right)^i,\chi)= \chi(-1)p^{m-1}J(\chi \left({ *\over p}\right)^{i}, \left({ *\over p}\right)) J_{m-2}(\chi \left({ *\over p}\right)^i,\chi).$$
}

\begin{proof}
Let $m$ be odd. 
 Then, by (1) of Proposition 5.8, we have
$$J_m(\chi\left({ *\over p}\right)^i ,\chi)=\left({-1 \over p}\right)^{(m-1)/2}p^{(m-1)/2}$$
$$\times \{J(\chi\left({ *\over p}\right)^i  ,\chi^{m})J_{m-1}(\chi\left({ *\over p}\right)^{i+1},\chi)+\chi(-1) I_{m-1}(\chi\left({ *\over p}\right)^{i+1},\chi)\}.$$
Thus the assertion holds if $\chi^m \not=1.$ Assume that $\chi^m=1.$ Then  by (2) of Proposition 5.8 and (2) of Proposition 5.7 we have  
$$J_{m-1}(\chi\left({ *\over p}\right)^{i+1},\chi)= \left({-1 \over p}\right)^{(m-1)/2}p^{(m-3)/2} J(\chi\left({ *\over p}\right)^{i},\left({ *\over p}\right)) $$
$$\times J(\chi\left({ *\over p}\right)^{i} ,\chi^{m-1} \left({ *\over p}\right)^{i})J_{m-2}(\chi \left({ *\over p}\right)^{i} ,\chi).$$
and 
$$I_{m-1}(\chi\left({ *\over p}\right)^{i+1},\chi)=\left({-1 \over p}\right)^{(m-3)/2}  p^{(m-3)/2}(p-1)\chi(-1) \left({-1 \over p}\right)^{i+1}$$
$$ \times  J(\chi\left({ *\over p}\right)^{i},\left({ *\over p}\right)) J_{m-2}(\chi \left({ *\over p}\right)^{i} ,\chi).$$
We note that  $J(\chi\left({ *\over p}\right)^i  ,\chi^{m})=-1,\chi(-1)=1$ and 
$$J(\chi\left({ *\over p}\right)^{i} ,\chi^{m-1} \left({ *\over p}\right)^{i})=J(\chi\left({ *\over p}\right)^{i} ,\overline{\chi \left({ *\over p}\right)^{i}})=-\chi(-1) \left({ -1   \over p}\right)^{i}=\left({ -1   \over p}\right)^{i+1}.$$
This proves the assertion.

Let $m$ be even. 
Then, by (2) of Proposition 5.8, we have
$$J_m(\chi\left({ *\over p}\right)^i ,\chi)=\left({-1 \over p}\right)^{(m-2)/2}p^{(m-2)/2}J(\chi\left({ *\over p}\right)^i ,\left({ *\over p}\right))  $$
$$\times \{J(\chi\left({ *\over p}\right)^i ,\chi^{m}\left({ *\over p}\right)^{i +1})J_{m-1}(\chi\left({ *\over p}\right)^{i+1},\chi)+\chi(-1) I_{m-1}(\chi\left({ *\over p}\right)^{i+1},\chi)\}.$$
Thus the assertion holds if $\chi^{m}\left({ *\over p}\right)^{i +1} \not=1.$ Assume that $\chi^{m}\left({ *\over p}\right)^{i +1} =1.$  Then by (1) of Proposition 5.7 and (1) of Proposition 5.8, we have
$$J_{m-1}(\chi\left({ *\over p}\right)^{i+1} ,\chi)=\left({-1 \over p}\right)^{(m-2)/2}p^{(m-2)/2}$$
$$ \times J(\chi \left({ *\over p}\right)^{i+1},\chi^{m-1})J_{m-2}(\chi \left({ *\over p}\right)^{i} ,\chi),$$
and 
$$I_{m-1}(\chi\left({ *\over p}\right)^{i+1} ,\chi)=\left({-1 \over p}\right)^{(m-2)/2}p^{(m-2)/2} J_{m-2}(\chi \left({ *\over p}\right)^{i} ,\chi)(p-1).$$
We note that  $J(\chi\left({ *\over p}\right)^i  ,\chi^{m}\left({ *\over p}\right)^{i +1})=-1,\left({ -1   \over p}\right)^{i+1}=1$ and 
$$J(\chi\left({ *\over p}\right)^{i+1} ,\chi^{m-1})=J(\chi\left({ *\over p}\right)^{i+1} ,\overline{\chi \left({ *\over p}\right)^{i+1}})=-\chi(-1) \left({ -1   \over p}\right)^{i+1}=-\chi(-1).$$
This proves the assertion.

\end{proof}

\bigskip

\noindent
{\bf Corollary.} {\it
 Let $\chi$ be a primitive character with an odd square free conductor $N.$ Assume that $\chi^2$ is primitive.  Then the value  $J_m(\chi)$ is nonzero. 
}
\begin{proof}
The assertion follows directly from the above theorem if $N$ is an odd prime. In general case, the assertion can also be proved by remarking that
$J_m(\chi)= \prod_{p | N} J_m(\chi^{(p)})$ and that ${\chi^{(p)}}^2 \not=1$ for any $p |N.$ 
\end{proof}
To compare our present result with the result in \cite{K-M}, we give the following:

\bigskip

\noindent
{\bf Proposition 5.10.}
 {\it  Let $\chi$  be  a primitive Dirichlet character mod $p.$ Assume that $\chi^2 \not=1.$ Then we have
$$J(\chi, \left({ * \over p}\right))J(\chi \left({ * \over p}\right), \chi \left({ * \over p}\right))=\left({ -1 \over p}\right) \bar \chi(4)p.$$ 
}
\begin{proof}
The assertion follows from [\cite{B-E}, Theorems 2.3 and 2.4].
\end{proof}

\noindent
{\bf Remark.} The above proof is due to the referee.  We note that the asserion can also be proved by using the same argument as in the proof of Theorem 5.5, and Lemma 5.3. 

By virtue of the above proposition, we see that Theorem 5.6 coincides with  \cite{K-M}, Proposition 3.7 in case $m=2.$
 
 Now let 
 $$F(Z)= \sum_{A \in {\mathcal L}_{n\ge 0}} c_F(A) {\bf e}({\rm tr}(AZ))$$
be an element of ${\textfrak M}_k(Sp_n({\bf Z}))$ and let $\chi$ be a Dirichlet character mod $N.$ Assume  $N \not= 2.$ Then by  [\cite {K-M}, Proposition 3.1], we have
$$L(s,F,\chi)=\sum_{A \in {\mathcal L}_{n>0}/SL_n({\bf Z})} {c_F(A)h(A,\chi) \over e(A)(\det A)^s}.$$
Thus by Theorem 5.6  we easily obtain:

\bigskip

\noindent
{\bf Theorem 5.11.} {\it
Let $N,p_i,l_i,u_{0,i} \ (i=1,\cdots,r)$ and $\chi$ be as in Theorem 5.6, and let $F$ be an element of ${\textfrak M}_k(Sp_n({\bf Z})).$ 

\noindent
{\rm (1).} If $\chi^{(p_i)}(u_{0,i}) \not=1$ for some $i.$ Then we have $L(s,F,\chi)=0.$ 

\noindent
{\rm (2).} Assume that $\chi^{(p_i)}(u_{0,i})=1$ for any $i.$ Fix a character $\tilde \chi$ such that  $\tilde \chi^n=\chi.$ 

\noindent
{\rm (2.1)} Let $n$ be even. Then we have 
$$L(s,F,\chi)= \prod_{i=1}^r (-1)^{(n-2)(p_i-1)/4 } \gamma_{n,p_i} $$
$$ \times \sum_{ \eta \in {\mathcal M}_{N,n} }  \overline {(\widetilde \chi \eta)}(2^n) \overline{J(\widetilde \chi \eta,\left({* \over N} \right))}\overline{J_{n-1}(\widetilde \chi \eta) }  L^*(s,F,\widetilde \chi \eta) .$$
{\rm (2.2)} Let $n$ be odd, and assume that $\chi^2 \not=1.$  Then we have 
$$L(s,F,\chi)= \prod_{i=1}^r (-1)^{(n-1)(p_i-1)/4 }\gamma_{n,p_i}  $$
$$ \times \sum_{ \eta \in {\mathcal D}_{N,m}} \overline {(\widetilde \chi \eta)}(2^{n-1}) \overline{J_{n-1}(\widetilde  \chi \eta)} L^*(s,F, \widetilde \chi \eta).$$
}

\section{Twisted Koecher-Maa{\ss} series of the first kind of  the D-I-I lift}

By Theorems 4.1 and 5.11, we obtain the following.

\bigskip

\noindent
{\bf Theorem 6.1.} {\it
Let $k$ and $n$ be positive even integers such that $n \ge 4, \ 2k-n \ge 12.$   
Let $h(z)$ and $E_{n/2+1/2}$ be as in Section 4. Let $N$ be a square free odd integer, and $N=p_1 \cdots p_r$ be the prime decomposition of $N.$ For each $i=1,\cdots,r$ let 
$l_i={\rm GCD}(n,p_i-1)$ and $u_{i} \in {\bf Z}$ be a  primitive $l_i$-th root of unity mod $p_i.$\\
 {\rm (1)}  Assume $\chi^{(p_i)}(u_i) \not=1$ for some $i.$ Then $L(s,I_n(h),\chi)=0.$ \\
{\rm (2)}  Assume $\chi^{(p_i)}(u_i) =1$ for all $i.$ Then 
$$L(s,I_n(h),\chi)=2^{ns} \sum_{\eta \in {\mathcal D}_{N,n}} \overline {\widetilde \chi \eta}(2^{n})\overline{J(\widetilde \chi \eta,\left({* \over N}\right) )} \overline{J_{n-1}(\widetilde \chi \eta)}$$
$$\times  \{ c_{n,N }R(s,h,E_{n/2+1/2},\widetilde \chi \eta)\prod_{j=1}^{n/2-1} L(2s-2j,S(h),(\widetilde \chi \eta)^2) $$
$$+d_{n,N}c_h(1)   \prod_{j=1}^{n/2} L(2s-2j+1,S(h),(\widetilde \chi \eta)^2)\}, $$
where $c_{n,N}$ and $d_{n,N}$ are nonzero rational numbers depending only on $n$ and $N,$ and $\widetilde \chi$ is a character such that $\widetilde \chi^n=\chi.$
} 
\bigskip

\noindent
{\bf Remark.}  In the case $n=2$, an explicit formula for $L(s,I_2(h),\chi)$ was given by  Katsurada-Mizuno \cite {K-M}.

\section{Applications}

Let $h_1$ and $h_2$ be modular forms of weight $k_1+1/2$ and $k_2+1/2,$ respectively, and $\chi$ be a Dirichlet character
. In Section 2, we reviewed on the algebraicity of the values $\widetilde R(m,h_1,h_2,\chi)$ at  half integers. We then naturally ask the following question:

\noindent
{\bf Question.} {What can one say about the algebraicity of  $\widetilde R(m,h_1,h_2,\chi)$ with $m$ an integer?}

As an application of Theorem 6.1, we give a partial answer to this question. We note that  
$$R(s,h_1,h_2,\chi)=(1-2^{-2s+k_1+k_2-1}\chi^2(2))^{-1}\widetilde R(s,h_1,h_2,\chi)$$
if the conductor of $\chi$ is odd.  Hence it suffices to consider the above question for $R(m,h_1,h_2,\chi)$ with integer $m$
 if $k_1+k_2$ is even.
 
Let $k$ and $n$ be positive even integers such that $n \ge 4, \ 2k-n \ge 12.$   
Let $h(z)$ and $E_{n/2+1/2}$ be as in Section 4. For a Dirichlet character $\chi$ of odd square free conductor $N=p_1\cdots p_r,$ we define 
$$R^{(\chi)}(s,h,E_{n/2+1/2})= \sum_{ \eta \in {\mathcal D}_{N,n}} \overline {\chi \eta}(2^{n})\overline{J(\chi \eta,\left({* \over N}\right)) } \overline{J_{n-1}(\chi \eta)}$$
$$\times R(s,h,E_{n/2+1/2},\chi \eta) \prod_{j=1}^{n/2-1} L(2s-2j,S(h),(\chi \eta)^2).$$

\noindent
{\bf Theorem 7.1.} {\it 
There exists a finite dimensional  $\overline {\bf Q}$-vector space 
$W_{h,E_{n/2+1/2}}$ in ${\bf C}$ such that
$$ {R^{(\chi)}(m,h,E_{n/2+1/2}) \over \pi^{mn}} \in W_{h,E_{n/2+1/2}}$$
for any integer $n/2+1 \le m \le k-n/2-1$ and all characters $\chi$ of odd square free conductor such that $\chi^n$ is primitive.
}
\begin{proof}
Put 
$${\bf M}^{(\chi)}(s,S(h))= \sum_{ \eta \in {\mathcal D}_{N,n}} \overline{ \chi \eta}(2^n) \overline{J(\widetilde \chi \eta,\left({* \over N}\right) )}\overline{J_{n-1}( \chi \eta)}$$
$$ \times \prod_{j=1}^{n/2} L(2s-2j+1,S(h),(\chi \eta)^2).$$
Then by Corollary to Proposition 3.1, we have  
 $${{\bf M}^{(\chi)}(m,S(h)) \over \pi^{mn}} \in  \overline {\bf Q}u_{-}(S(h))^{n/2} \pi^{-n^2/4}.$$
 By Theorem 6.1, we have 
 $$ L(m,I_n(h),\chi^n)$$
$$= 2^{nm}\{ c_{n,N}   R^{(\chi)}(m,h,E_{n/2+1/2}) +d_{n,N}c_h(1){\bf M}^{(\chi)}(m,S(h))\}.$$
Hence by Theorem 2.2, we have 
 $${R^{(\chi)}(m,h,E_{n/2+1/2}) \over \pi^{mn}}  \in \overline {\bf Q}u_1  \otimes_{\overline {\bf Q}} V_{I_n(h)}+\overline {\bf Q}u_2$$ 
with some complex numbers $u_1$ and $u_2,$ where $V_{I_n(h)}$ is the $\overline {{\bf Q}}$-vector space associated with
 $I_n(h)$ in Theorem 2.2. This proves the assertion.
\end{proof}

By the above theorem, we immediately obtain the following:

\noindent
 {\bf Theorem 7.2.} {\it
 Let $d > \dim_{\overline {\bf Q}} W_{h,E_{n/2+1/2}}.$ Let $m_1,m_2,\cdots,m_d$ be integers such that $n/2+1 \le m_1,m_2,\cdots, m_d  \le k-n/2-1$ and $\chi_1,\chi_2,\cdots,\chi_d$ be Dirichlet characters of odd square free  conductors $N_1,N_2,\cdots,N_d,$ respectively such that $\chi_i^n$ is primitive for any $i=1,2,\cdots d.$ Then the values ${\displaystyle R^{(\chi_1)}(m_1,h,E_{n/2+1/2}) \over \displaystyle \pi^{m_1 n}} ,\cdots,  $\\
${\displaystyle R^{(\chi_d)}(m_d,h,E_{n/2+1/2}) \over \displaystyle \pi^{m_d n}}$ are linearly dependent over $\overline {\bf Q}.$
}
\bigskip

\noindent
 {\bf Corollary.} {\it 
 In addition to the notation and the assumption as above, assume that $n \equiv 2 \ {\rm mod} \ 4.$ Then the values $$\left\{{\displaystyle R(m_i,h,E_{n/2+1/2}, \chi_i \eta_{ij})  \over \displaystyle\pi^{2m_i}} \right \} _{1 \le i \le d,  \eta_{i j} \in {\mathcal D}_{N_i,n}} $$ 
 are linearly dependent over $\overline{ {\bf Q}}.$ 
}
 \begin{proof} Put ${\bf L}_n(s,S(h),(\chi_i \eta_{ij})^2)=\prod_{l=1}^{n/2-1} L(2s-2l,S(h),(\chi_i \eta_{ij})^2).$ Then by Theorem 1.1, the value ${{\bf L}_n(m_i,S(h),(\chi_i \eta_{ij})^2)  \over \pi^{m_i(n-2)}}$ belongs to $\overline {{\bf Q}}u_{+}(S(h))^{n/2-1}\pi^{-n^2/4+n/2},$
 and in particular if $n \equiv 2 \ {\rm mod} \ 4, $ then  it is nonzero for any $\chi.$ Moreover, by Corollary to Theorem 5.10, $J(\chi_i \eta_{ij},({* \over N}))J_{n-1}(\chi_i \eta_{ij})$ is non-zero and belongs to $\overline {{\bf Q}}.$ 
 Thus the assertion holds.
 \end{proof}

As another application of Theorem 6.1, we also have a functional equation for $R^{(\chi)}(s,h,E_{n/2+1/2}).$ Namely, by Theorem 2.1 we obtain:
 
 \noindent
 { \bf Theorem 7.3.} {\it 
Let $h$ be as above. Let $\chi$ be a primitive character of odd square free conductor $N.$ Assume that $n \equiv 2 \ {\rm mod} \ 4,$ and that $\chi^n$ is primitive. Put 
 $${\mathcal R}^{(\chi)}(s,h,E_{n/2+1/2})=N^{2s}\tau(\chi^n)^{-1}\gamma_n(s)R^{(\chi)}(s,h,E_{n/2+1/2}).$$
 Then ${\mathcal R}^{(\chi)}(s,h,E_{n/2+1/2})$ has an analytic continuation to the whole $s$-plane, and has the following functional equation:
 $${\mathcal R}^{(\chi)}(k-s,h,E_{n/2+1/2})={\mathcal R}^{(\chi)}(s,h,E_{n/2+1/2}).$$
}
 
 \noindent
 {\bf Remark.} (1) As functions of $s,$ the Dirichlet series \\
$\left\{ R(s,h,E_{n/2+1/2}, \chi_i \eta_{ij}) \right \} _{1 \le i \le d,  \eta_{ij} \in {\mathcal D}_{N_i,n}} $
are linearly independent over ${\bf C}.$

 (2) In the case of $n=2,$ this type of result was given for $R(m,h,E_{3/2})$ with $E_{3/2}$ Zagier's Eisenstein series of weight $3/2$ by \cite{K-M}.
 
 (3) The meromorphy of this type of series was derived in \cite{Sh3} by using so called the Rankin-Selberg integral expression in more general setting, but we don't know whether  the functional equation of the above type can be directly proved without using the above method. 

{\bf Acknowledgement}
 
 The author thanks the referee for giving useful comments especially on Theorem 5.5 and Proposition 5.10, which make our paper consice.

\noindent
\author{Hidenori KATSURADA \\[1mm]
Muroran Institute of Technology \\
 27-1 Mizumoto, Muroran, 050-8585, Japan \\
E-mail: \verb+hidenori@mmm.muroran-it.ac.jp+ \\[2.5mm]}

\end{document}